\documentclass[11pt]{article}

\usepackage{amsfonts,amsmath,graphicx}
\graphicspath{{.}}
   \makeatletter
   \let\accent@spacefactor\relax
   \makeatother
   \usepackage[french]{babel}
\usepackage{amsfonts}
\usepackage{amsmath}   
\usepackage{amsthm}    
\usepackage{amscd}     
\usepackage{euscript}

\def\C{{\Bbb C}}
\def\P{{\Bbb P}}
\def\F{{\Bbb F}}

\def\Z{{\Bbb Z}}

\newtheorem{defi}{D\'{e}finition}[section]
\newtheorem{prop}[defi]{Proposition}

\newtheorem{lem}[defi]{Lemme}

\newtheorem{theo}[defi]{Th\'{e}or\`{e}me}

\newtheorem{coro}[defi]{Corollaire}

\newtheorem{rem}[defi]{Remarque}

\title{
 Fibr\'es vectoriels de rang deux sur $\P^2$\\
provenant d'un rev\^etement double}
\author{ Jean Vall\`es}
\date{19/06/08}

\begin{document}

\maketitle

\begin{abstract}
Depuis Schwarzenberger et son c\'el\`ebre article intitul\'e ``Vector bundles on the projective plane'', on sait que tout fibr\'e de rang deux sur $\P^2(\C)$ peut-\^etre d\'efini comme l'image directe d'un faisceau inversible sur une surface recouvrant doublement le plan. Ce th\'eor\`eme sugg\`ere d'\'etudier les fibr\'es de rang deux en fonction de la courbe de ramification du rev\^etement dont ils proviennent. 


Ainsi, dans la premi\`ere partie on d\'emontre que, \'etant donn\'e un rev\^etement 
ramifi\'e le long d'une courbe irr\'eductible  $C_{2r}$  
de degr\'e $2r$,  les  droites  sauteuses 
d'ordre fix\'e (ordre d\'ependant de $r$ et de la parit\'e du fibr\'e de rang deux) 
des fibr\'es images directes sont n\'ecessairement  $r$-tangentes \`a la courbe $C_{2r}$.


Dans la seconde partie nous nous penchons plus particuli\`erement sur le cas $r=2$. Nous donnons alors une liste de fibr\'es pour lesquels les droites sauteuses sont exactement les bitangentes de la quartique de ramification.


\begin{center}
\textbf{Abstract}
\end{center}
\hspace{0.3cm} Since Schwarzenberger and his celebrated paper called ``Vector bundles on the projective plane'' we know that any rank two vector bundle on $\P^2$ is a direct image of a line bundle on a double covering of the plane. This theorem suggests to study the rank two vector bundles according to the branch curve of the covering  which they come from.


Thus, in the first part we prove that, given a double covering ramified over an irreducible curve $C_{2r}$ with degree $2r$,  the jumping lines of fixed order (order depending on $r$ and on the parity of the rank two vector bundle) of the direct images vector bundles are necessarely  $r$-tangent to $C_{2r}$.


In the second part we concentrate on the case $r=2$. Then we give  a list of vector bundles for which the jumping lines are exactly the bitangent lines to the branch quartic.
\end{abstract}
{\small Mots cl\'es : fibr\'es de rang deux, droites de saut,  courbe de ramification, bitangentes}
\section{Introduction}

Tout fibr\'e de rang deux sur $\P^2(\C)$ peut-\^etre d\'efini comme l'image directe d'un faisceau inversible sur une surface recouvrant doublement le plan (\cite{S}, thm 3). Ce th\'eor\`eme est le point de d\'epart de cet article et son fil rouge. Reprenant donc l'id\'ee de Schwarzenberger, on se propose  d'\'etudier les fibr\'es de rang deux comme images de faisceaux inversibles d'un rev\^etement donn\'e. Plus pr\'ecis\'ement ce dernier \'etant caract\'eris\'e par sa ramification, on \'etudiera les images directes pour une courbe de ramification fix\'ee.

\smallskip

Consid\'erons donc   $\pi : X_r \rightarrow \P^2$ un rev\^etement double  du plan projectif complexe ramifi\'e le long d'une courbe lisse  $C_{2r}$ de degr\'e $2r$, $L$ un fibr\'e inversible sur $X_r$ et $E=\pi_{*}L$ son image directe dans $\P^2$.  Schwarzenberger montre, pour les cas $r=1$ et $r=2$, les deux propositions suivantes (on d\'efinit les notions de droites de saut et d'ordre de saut dans la section qui suit cette introduction)

\smallskip

\textbf{Proposition}(\cite{S}, prop.8, $r=1$) \textit{Une droite est une droite de saut de $E$ si et seulement si elle est tangente \`a la conique de ramification, que  $c_1(E) $ soit paire ou impaire.}

\smallskip

\textbf{Proposition}(\cite{S}, prop.10, $r=2$) \textit{Lorsque $c_1(E) $ est impaire toute droite de saut de $E$
est tangente \`a la quartique de ramification. Lorsque $c_1(E) $ est paire, ce n'est pas n\'ecessairement le cas.}

\smallskip
Une vingtaine d'ann\'ees plus tard Ottaviani g\'en\'eralisera ce r\'esultat en reprenant l'id\'ee de la d\'emonstration de la prop. 10 de \cite{S},

\smallskip

\textbf{Proposition}(\cite{O}, prop.6, $r\ge 1$) \textit{ 
\begin{center}
$l$ est une droite $(\lceil \frac{r+\epsilon}{2}\rceil+1)$-sauteuse de $E$ $\Longrightarrow $ $l$ est tangente \`a $C_{2r}$.
\end{center}
o\`u $\epsilon=0$ si $c_1(E)$ est paire et vaut $-1$ sinon.}

\smallskip

Si l'on \'evoque les tangentes d'une quartique lisse, on pense  obligatoirement \`a ses  $28$ bitangentes (celles-ci  d\'eterminent la quartique, d'apr\`es \cite{CS}, qui \`a son tour d\'etermine le rev\^etement). Quelle ne fut pas ma surprise en d\'ecouvrant que, dans le cas $r=2$ et lorsque $c_1(E)$ est impaire, toute  droite de saut de $E$ est non seulement tangente mais m\^eme  bitangente \`a la quartique, concentrant le sch\'ema des droites de saut sur un ensemble fini! Plus  g\'en\'eralement, on montre  

\smallskip

\textbf{Th\'eor\`eme} (thm \ref{bt}, section \ref{Ip})\textit{
\begin{center}
$l$ est une droite $(\lceil \frac{r+\epsilon}{2}\rceil+1)$-sauteuse de $E$ $\Longrightarrow $ $l$ est $r$-tangente \`a $C_{2r}$.
\end{center}
o\`u $\epsilon=0$ si $c_1(E)$ est paire et vaut $-1$ sinon.}

\smallskip

\textbf{Remarque :} De plus, pour prouver ce dernier r\'esultat, l'hypoth\`ese d'irr\'eductibilit\'e de la courbe de ramification (au lieu de la lissit\'e) suffit.

\smallskip

Une cons\'equence int\'eressante de ce th\'eor\`eme est que  pour $r=3$, lorsque la sextique de ramification est g\'en\'erale, en particulier lorsqu'elle ne  poss\`ede pas de tritangente, les fibr\'es images directes n'ont pas de droites bisauteuses. Ce th\'eor\`eme permet aussi  de dresser une liste des fibr\'es  instables (corollaire \ref{ins}) ou d\'ecompos\'es (corollaire \ref{dec}) pouvant s'\'ecrire comme images directes. La preuve, d'inspiration diff\'erente de celle de Schwarzenberger, exploite l'existence de la  section non nulle d'un d\'ecal\'e de  $S^2E$ qui donne le plongement $X_r\hookrightarrow \P E$ (lemme  \ref{lece}).

\smallskip

Dans une deuxi\`eme partie (cf section \ref{IIp}),  d\'evolue plus particuli\`erement au cas $r=2$, on donne  une liste de fibr\'es stables  de premi\`ere classe de Chern impaire dont l'ensemble des droites de saut est exactement  l'ensemble des bitangentes de la quartique (th\'eor\`eme \ref{qpon} et la remarque \ref{equivalence}).  Pour d\'ecrire les sch\'emas de droites de saut, notamment les ordres de saut, nous introduisons des suites exactes courtes naturelles reliant les images directes entre elles. Cette seconde partie, notamment le calcul explicite des images directes,  s'appuie sur 
une description d\'etaill\'ee de l'involution bien connue de Geiser (plus exactement il s'agit de la transformation appel\'ee $(1,2)$ par Bateman dans \cite{B}) : le rev\^etement double est construit \`a partir de $7$ points en position g\'en\'erale du plan et du syst\`eme lin\'eaire des cubiques ayant ces points comme lieu de base. 

\medskip

{\small Je  remercie chaleureusement les membres de l'\'equipe de g\'eom\'etrie alg\'ebrique de l'universit\'e Ulisse Dini de Florence qui m'ont accueilli et tout particuli\`erement Elena Rubei pour son invitation et Giorgio Ottaviani pour son soutien permanent.
Si j'ai eu le loisir de travailler \`a Florence et d'y \'ecrire ce texte en ayant le sentiment que (...)\textit{``tout ce que je fis durant mon s\'ejour ne fut en effet que l'occupation d\'elicieuse et
n\'ecessaire d'un homme qui s'est d\'evou\'e \`a l'oisivet\'e''} \footnote{Les r\^everies du promeneur solitaire (J.J. Rousseau)} c'est aussi gr\^{a}ce \`a l'obtention d'une d\'el\'egation au CNRS dont je tiens \`a souligner combien elle est bienvenue dans la carri\`ere d'un ma\^itre de conf\'erences.}

\section{Rappels g\'en\'eraux et notations principales}

1) \textbf{Droites de saut.}  Soit $E$ un fibr\'e vectoriel de rang deux sur le plan projectif complexe $\P^{2}$. Sur toute droite 
$l\in \P^{2\vee}$ le fibr\'e restreint $E_l$ est une somme directe de faisceaux inversibles : il existe 
$a_l,b_l\in \Z$ v\'erifant $$ E_l=O_l(a_l)\oplus O_l(b_l), \, \, a_l+b_l=c_1(E)$$
Afin de limiter les tournures un peu lourdes on appellera parfois \textit{fibr\'es pairs} (respectivement \textit{impairs}) les fibr\'es vectoriels de rang deux dont la premi\`ere classe de Chern est paire (respectivement impaire). On associe au fibr\'e $E$ l'entier $\epsilon:=\epsilon(E)$ qui vaut $0$ lorsque $E$ est pair et $-1$ lorsque $E$ est impair. 
\begin{defi} Soit $i$ un entier strictement positif.
Une droite $l$ est une droite $i$-sauteuse de $E$ si et seulement si $\vert a_l-b_l\vert \ge 2i-\epsilon$.  
On note $S_i(E)$ l'ensemble de ses droites  $i$-sauteuses et $S(E)$ le sch\'ema de ses droites de saut. 
\end{defi}
Habituellement on r\'eserve cette appelation aux fibr\'es vectoriels de rang deux stables ou semi-stables, c'est \`a dire essentiellement (on  excepte $O_{\P^2}^2$) les fibr\'es vectoriels de rang deux  pour lesquels l'implication suivante est v\'erifi\'ee :
$$ c_1(E(n))<0 \Rightarrow h^0(E(n))=0$$
Dans ce cas, conform\'ement au th\'eor\`eme de Grauert-Mulich, la droite g\'en\'erale n'est pas sauteuse. De plus lorsque $E$
est pair  $S(E)$ est un diviseur (une courbe de $\P^{2\vee}$)
et lorsque $E$ est impair le sch\'ema est en g\'en\'eral fini 
mais peut aussi dans des cas particuliers contenir  une courbe; 
lorsqu'il est fini sa longueur vaut $\binom{c_2(E)}{2}$.

\smallskip

Cependant, tout comme dans l'article originel de Schwarzenberger, c'est la d\'ecomposition du  fibr\'e $E$ en fonction du fibr\'e en droite dont il provient qui nous int\'eresse, pas sa stabilit\'e. C'est pourquoi garder la m\^eme terminologie pour tous semblait moins d\'eroutant que d'en introduire une diff\'erente (elles sont appel\'es exceptionnelles par Scwharzenberger). Ainsi,  
pour les fibr\'es instables on aura (voir ci-dessous le rappel) $S_1(E)=\P^{2\vee}$ et $S_i(E)$ pour  un $i>1$ sera support\'e par une r\'eunion de droites de $\P^{2\vee}$. Pour les fibr\'es d\'ecompos\'es en somme directe de fibr\'es inversibles sur $\P^2$ le saut est constant.

\medskip

2)\textbf{Fibr\'es instables non d\'ecompos\'es.} Rappelons maintenant quelques r\'esultats classiques  concernant les fibr\'es instables sur le plan projectif. 

\smallskip

Soit $E$ un fibr\'e de rang deux sur $\P^2$, non d\'ecompos\'e comme somme directe de fibr\'es en droites et tel que   que $c_1(E)=-1$ ou $-2$. Alors $E$ est instable si et seulement s'il existe $k\ge 0$ tel que $H^0(E(-k))\neq 0 $ et $H^ 0(E(-1-k))=0$. De plus $h^0(E(-k))=1$ et cette unique section s'annule en codimension $2$ le long d'un sch\'ema non vide (sinon le fibr\'e est d\'ecompos\'e) de longueur $c_2(E(-k))$. Ainsi l'on peut \'ecrire   
$$ \begin{CD}
0 @>>> O_{\P^{2}} @>>>E(-k)
@>>> J_{Z}(-2k+c_1)
@>>>0
\end{CD}
$$
On constate que sur la droite g\'en\'erale $E_l=O_l\oplus O_l(c_1-2k)$ tandis que sur les droites rencontrant $Z$ en degr\'e $i$ on a $E_l=O_l(i)\oplus O_l(c_1-2k-i)$.

\section{Puissance sym\'etrique associ\'ee au rev\^etement}\label{Ip}
Commen\c{c}ons par un r\'esultat presque tautologique mais tr\`es utile.
\begin{lem}\label{lece}
Soit  $\pi : X_r\rightarrow \P^2$ un rev\^etement double  ramifi\'e le long d'une courbe irr\'eductible de degr\'e $2r$ de $\P^2$. Soient $L$ un faisceau inversible sur la surface $X_r$ et  $E=\pi_*{L}$ fibr\'e de rang deux sur $\P^2$. Alors il existe une section non nulle $s\in H^0((S^2E)(r-c_1(E)))$ et plus pr\'ecis\'ement  un homomorphisme sym\'etrique $\phi_s \in Hom(E^{\vee},E(r-c_1(E)))$ dont le conoyau est un faisceau de rang $1$ support\'e par la courbe de ramification
\end{lem}
\begin{proof}
On note
$h$ la classe de $O_{{\P}^2}(1)$, $\eta $  la classe du relativement ample
sur le fibr\'e projectif
$\P E$,  $\pi$ le morphisme canonique
$ \pi :\P E \longrightarrow \P^2$ et $c_1$ la premi\`ere classe de Chern de $E$.

\smallskip

 Une section non nulle  $s\in
H^0(S^2E(r-c_1))=H^0(O_{\P E}(2\eta+(r-c_1)h))$  d\'efinit un diviseur   $X_r$ de $\P E$  qui
est un rev\^etement double de
$\P^2$. En effet consid\'erons  la suite exacte
d\'efinissant le diviseur $X_r$  
$$ 0\rightarrow O_{\P E}(-2\eta+(c_1-r)h)
\longrightarrow O_{\P E}\longrightarrow O_{X_r} \rightarrow
0$$
Le faisceau dualisant relatif est $\omega_{\P E/\P^2}=O_{\P E}(-2\eta+c_1h)$ (\cite{Ha}, ex. 8.4). On en d\'eduit que 
$R^1\pi_*O_{\P E}(-\eta)=0$ et  $R^1\pi_*O_{\P E}(-2\eta)=
O_{\P^2}(-c_1)$, ce qui implique  $$E=\pi_*O_{X_r}(\eta) \,\,\mathrm{et} \,\, \pi_*O_{X_r}=O_{\P^2}\oplus O_{\P^2}(-r)$$ 
R\'eciproquement, 
notons $f=0$ l'\'equation de la courbe de ramification. Alors, dans un sens plus imag\'e que  pr\'ecis (pour plus de pr\'ecision on consultera, par exemple, les pages 173-174 de \cite{Pe}), on peut consid\'erer que   la racine carr\'ee $\sqrt{f}$ existe en tant que diviseur $D$ sur $X_r$ et qu'elle induit une suite exacte, canonique (o\`u $D$ est la courbe de ramification ``en haut'')
$$ 0\rightarrow  \pi ^{*}O_{\P^2}(-r)
\longrightarrow O_{X_r}\longrightarrow O_D \rightarrow
0$$
Cette derni\`ere suite exacte, tensoris\'ee par un fibr\'e en droite $L$ de $X_r$, descend sur $\P^2$ 
\begin{equation}
\label{eq}
0\rightarrow \pi_*L(-r)
\longrightarrow \pi_*L\longrightarrow K \rightarrow 0
\end{equation}
o\`u $K$ est un faisceau de rang 1 support\'e par la courbe de ramification $\{ f=0\}$. Par ailleurs cette application est sym\'etrique, puisque  la ramification est le lieu des points fixes sous l'action du groupe de Galois $\Z/2\Z$. Pour se convaincre de la sym\'etrie on utilisera la d\'ecomposition suivante des fibr\'es $$ E\otimes E=\bigwedge^2 E\oplus S^2E $$
o\`u $E=\pi_*L$.
La suite exacte (\ref{eq}) correspond \`a une section non nulle de $E\otimes E(r-c_1)$. 
La d\'ecomposition canonique en parties sym\'etrique et antisym\'etrique s'\'ecrit
$$ E\otimes E(r-c_1)=O_{\P^2}(r)\oplus S^2E(r-c_1) $$ Ceci donne sur les espaces de sections   
$$ H^0(E\otimes E(r-c_1))=H^0(O_{\P^2}(r))\oplus H^0((S^2E)(r-c_1)) $$
Soit $P$ une forme non nulle de degr\'e $r$. Elle induit une section $s_P\in H^0(E\otimes E(r-c_1))$ d\'ecrite ci-dessous
$$ \begin{CD}
0 @>>>E(-r) @>s_P>>E
@>>> E_{\lbrace P=0\rbrace}
@>>>0
\end{CD}
$$
Le conoyau est un faisceau de rang deux sur la courbe $\lbrace P=0 \rbrace$ de degr\'e $r$ et non pas un faisceau  de rang $1$ support\'e par une courbe de degr\'e $2r$.  Ceci montre que la section de $H^0(E\otimes E(r-c_1))$ obtenue via la racine carr\'ee de $f$ provient de la partie sym\'etrique $H^0((S^2E)(r-c_1))$. 
\end{proof}
Gr\^ace \`a ce rappel nous prouvons, 
\begin{theo}\label{bt} Soient $\pi : X_r \rightarrow \P^2$ un rev\^etement double  ramifi\'e le long d'une courbe irr\'eductible  $C_{2r}$ de degr\'e $2r$, $L$ un fibr\'e inversible sur $X_r$ et $E=\pi_{*}L$ son image directe dans $\P^2$.  Alors, 
\begin{center}
$l$ est une droite $(\lceil \frac{r+\epsilon}{2}\rceil+1)$-sauteuse de $E$ $\Longrightarrow $ $l$ est $r$-tangente \`a $C_{2r}$.
\end{center}
o\`u $\epsilon=0$ si $c_1(E)$ est pair et  vaut $-1$ sinon.
\end{theo} 
\begin{rem}
En d'autres termes, quitte \`a suppposer que $c_1(E)=0$ ou $-1$, cet \'enonc\'e affirme qu'une droite $l$ au dessus de laquelle :\\
  $E_l=O_l(s+1)\oplus O_l(-s-1)$ lorsque  $r=2s$ et $c_1(E)=0$ ou\\
  $E_l=O_l(s)\oplus O_l(-s-1)$ lorsque  $r=2s$  et $c_1(E)=-1$ ou\\
 $E_l=O_l(s+1)\oplus O_l(-s-1)$ lorsque  $r=2s+1$  et $c_1(E)=0$ ou\\
 $E_l=O_l(s+1)\oplus O_l(-s-2)$ lorsque  $r=2s+1$  et $c_1(E)=-1$ \\
est n\'ecessairement $r$-tangente \`a la courbe $C_{2r}$.
\end{rem}

\begin{proof}
Consid\'erons la section non nulle de $ H^0((S^2E)(r-c_1))$ donn\'ee par le lemme pr\'ec\'edent. Elle induit une application sym\'etrique $$ \phi_s: E^{\vee} \longrightarrow E(r-c_1)$$ 
Si $E_l=\mathcal{A}\oplus \mathcal{B}$ (avec $\mathcal{A}$ et $\mathcal{B}$ fibr\'es inversibles sur $l$) nous pouvons \'ecrire 
$$ (\phi_s)_{\mid l} : \mathcal{A}^{\vee} \oplus \mathcal{B}^{\vee} \rightarrow \mathcal{A}(r)\oplus \mathcal{B}(r)$$
ou bien 
$$ (\phi_s)_{\mid l} : \mathcal{B}^{\vee}\oplus \mathcal{A}^{\vee} \rightarrow \mathcal{A}(r)\oplus \mathcal{B}(r)$$

\smallskip

--- Montrons que le premier cas ne peut se produire. En effet, compte tenu de la sym\'etrie de $\phi_s$, il m\`ene \`a 
$$ \begin{CD}
E^{\vee}_l=O_l(-i)\oplus O_l(i-c_1) @> {\left (
\begin{array}{cc}
            E_l &0\\
             0  & F_l
\end{array}
  \right )}>> E_l(r-c_1)=O_l(r-i)\oplus O_l(r+i-c_1)
\end{CD}
$$
Avec 
$E_l, F_l$ des polyn\^omes de degr\'es \'egaux \`a $r$.

\smallskip

Soit alors, $f(x_0,x_1,x_2)=0$ l'\'equation de la courbe de ramification. Pla\c{c}ons-nous sur un ouvert de Zariski $U$, disons pour simplifier, sur l'ouvert
$U=\lbrace x_0\neq 0 \rbrace$. Sur cet ouvert, l'application $\phi : E\rightarrow E(r)$ est repr\'esent\'ee par une matrice $$M_U=  \left (
\begin{array}{cc}
            E &G\\
              G & F
\end{array}
  \right )$$ o\`u $E,F,G$ sont des polyn\^omes (non homog\`enes) en les variables affines 
$(x_1,x_2)$. Son  d\'eterminant vaut $EF-G^2=f(x_1,x_2)$, o\`u l'on note encore $f(x_1,x_2)=0$ l'\'equation affine de la courbe de ramification. 

\smallskip

On remarque que, except\'e  lorsque $E=O_{\P^2}^2$,  il existe toujours, pour $E$ semi-stable, stable, instable ou d\'ecompos\'e, au moins un diviseur de droites (i.e. une famille une dimensionnelle de droites) et un entier $i>0$ tels que 
$$ E_l=O_l(i+c_1)\oplus O_l(-i) $$
Par exemple, si $E$ est stable c'est vrai, sur la droite g\'en\'erale pour $i=1, c_1=-1$ et sur un diviseur pour $i=1, c_1=0$.

\smallskip

On consid\`ere donc la restriction aux  droites $l\subset U$ sur lesquelles 
$$E_l=O_l(i+c_1)\oplus O_l(-i), \,\, i>0$$ (elles forment un diviseur de $U$). On a 
$$ \begin{CD}
E_l=O_l(i+c_1)\oplus O_l(-i) @> {\left (
\begin{array}{cc}
            E_l &G_l\\
              G_l & F_l
\end{array}
  \right )}>> E_l(r)=O_l(r+i+c_1)\oplus O_l(r-i)
\end{CD}
$$
Mais alors  $\mathrm{deg}(G_l)=r-2i-c_1$ et $\mathrm{deg}(G_l)=r+2i+c_1$. 
Ceci implique $G_l\equiv 0$ sur la famille une-dimensionelle de droites de l'ouvert $U$, par cons\'equent
$G\equiv 0$ sur l'ouvert $U$. On en d\'eduit que $f(x_1,x_2)=0$ n'est pas irr\'eductible ce qui contredit l'hypoth\`ese du th\'eor\`eme.

\medskip

--- Ce premier cas s'av\'erant impossible, consid\'erons la seconde \'ecriture. L'application $\phi_s$ restreinte \`a la droite $l$, compte tenu de la sym\'etrie, est donc 
$$ \begin{CD}
E^{\vee}_l=O_l(i-c_1)\oplus O_l(-i) @> {\left (
\begin{array}{cc}
            A_l &B_l\\
              B_l & C_l
\end{array}
  \right )}>> E_l(r-c_1)=O_l(r-i)\oplus O_l(r+i-c_1)
\end{CD}
$$
Avec $A_l,B_l,C_l$ des polyn\^omes sur $l$ de degr\'es respectivement \'egaux \`a $r-2i+c_1$, $r$ et $r+2i-c_1$.

\smallskip

L'hypoth\`ese ``$l$ est une droite $(\lceil \frac{r+\epsilon}{2}\rceil+1)$-sauteuse de $E$'' implique que $r-2i+c_1<0$. Mais alors cette in\'egalit\'e  impose que $A_l\equiv 0$. Par cons\'equent le sch\'ema support\'e par l'intersection $l\cap C_{2r}$ est d\'efini par l'\'equation $B_l^2=0$. Et le th\'eor\`eme est prouv\'e.
\end{proof}
\begin{rem}
Ainsi dans le cas $r=1$ les droites de saut co\"{\i}ncident avec les droites tangentes \`a la conique $C_2$ que $E$ soit pair ou impair. Lorsque $r=2$ les droites de saut, pour $E$  impair, et les droites bisauteuses, pour $E$ pair, sont obligatoirement bitangentes \`a la quartique de ramification $C_4$ (et donc, ensemblistement, elles sont au plus $28$ droites distinctes!)  Enfin pour $r=3$, une sextique g\'en\'erale $C_3$ n'ayant pas de tritangente, les fibr\'es associ\'es, pairs ou impairs,  n'ont pas de bisauteuses.
\end{rem}
Ce th\'eor\`eme valant pour toutes les images directes de faisceaux inversibles sans consid\'eration de stabilit\'e, il permet de dresser une liste des fibr\'es instables et  d\'ecompos\'es provenant d'en haut.
\begin{coro}\label{dec}
Soit  $L$ un faisceau inversible  sur $X_r$ tel que $E=\pi_{*}L$ est d\'ecompos\'e. Alors  au d\'ecalage pr\`es on a $E= O_{\P^2}\oplus O_{\P^2}(-i)$ avec $0\le i\le r$.
\end{coro}
\begin{proof}
Si $i$ v\'erifiait $i>r$  toutes les droites de $\P^2$ devraient \^etre $r$-tangentes \`a la ramification, ce qui n'est pas.
\end{proof}
\begin{coro}\label{ins}
Soit  $L$ un faisceau inversible  sur $X_r$ tel que $E=\pi_{*}L$. On suppose qu'il n'est pas d\'ecompos\'e. Alors
$$ i>\lceil \frac{r}{2} \rceil+\epsilon \Rightarrow h^0(E(-i))=0 $$
\end{coro}
\begin{rem}
1) Pour $r\ge 3$, lorsque $C_{2r}$ n'a pas de tritangente, l'unique section non nulle de 
$h^0(E(-\lceil \frac{r}{2}\rceil-\epsilon))$ s'annule sur un point lisse, en d'autres termes les seuls cas, apr\`es normalisation, sont 
$$(c_1,c_2)=(-1,1-\frac{(r-3)(r-1)}{4})\,\, \mathrm{ou}\,\, (c_1,c_2)=(0,1-\frac{r^2}{4})$$

2) Pour $r=1$ et $r=2$ les seules images directes (de faisceaux inversibles) instables sont d\'ecompos\'ees.
\end{rem}
\begin{proof}
La preuve pour $c_1=-1$ et $r $ pair suffit \`a se convaincre de la validit\'e du corollaire.
Pla\c{c}ons nous directement sur le bord et consid\'erons l'unique section non nulle associ\'ee
$$ \begin{CD}
0 @>>> O_{\P^{2}} @>>>E(-\frac{r}{2}+1)
@>>> J_{Z}(-r+1)
@>>> 0
\end{CD}
$$
Pour une droite $l$ g\'en\'erale on a $E_l(-\frac{r}{2}+1)=O_l\oplus O_l(-r+1)$. Si une droite $l$ intersecte le sch\'ema des z\'eros $Z$ (n\'ec\'essairement non vide par hypoth\`ese) en degr\'e $i$, la restriction du fibr\'e est 
$E_l(-\frac{r}{2}+1)=O_l(i)\oplus O_l(-r+1-i)$. Si $i>1$ la droite $l$ est $r$-tangente \`a la courbe de ramification
en vertu du th\'eor\`eme pr\'ec\'edent. Par cons\'equent, s'il n'existe pas de $r$-tangente \`a la courbe, le sch\'ema $Z$ est support\'e par un point et  il est lisse (de longueur \'egale \`a $1$).
\end{proof}

\subsection{Th\'eta-caract\'eristiques sur $C_{2r}$}
Soit $E$ un fibr\'e de rang deux obtenu  comme image directe d'un faisceau inversible sur $X_r$. Consid\'erons l'homomorphisme sym\'etrique correspondant \`a la section non nulle $s\in H^0((S^2E)(r-c_1(E)))$. On note $\theta$ le faisceau de rang $1$ qui apparait comme conoyau de cet homomorphisme. \`A la suite exacte 
$$ \begin{CD}
0 @>>>E(-r) @>>>E
@>>> \theta
@>>>0
\end{CD}
$$
on applique le foncteur $\underline{Hom}(.,O_{\P^2}(-3))$ pour avoir
$$ \begin{CD}
0 @>>>E^{\vee}(-3) @>>>E^{\vee}(r-3)
@>>> \underline{Ext}^1(\theta,O_{\P^2}(-3))
@>>>0
\end{CD}
$$
ce qui prouve $\underline{Ext}^1(\theta,O_{\P^2}(-3))=\theta(r-c_1-3)=\theta^{*}\otimes \omega_{C_{2r}}$. Comme 
$ \omega_{C_{2r}}=O_{C_{2r}}(2r-3) $ on trouve 
$\theta^2=O_{C_{2r}}(r+c_1)$, i.e le faisceau conoyau est une racine carr\'ee de $O_{C_{2r}}(r+c_1)$. Quand $c_1=r-3$  c'est une th\'eta-caract\'eristique sur $C_{2r}$. On la retrouvera dans la deuxi\`eme partie.

\section{En famille, autour d'une quartique plane}\label{IIp}
Dans cette partie nous nous concentrons sur le cas $r=2$ dont Schwarzenberger a commenc\'e l'\'etude dans \cite{S}.
On rappelle tout d'abord la construction du rev\^etement double comme \'eclatement de sept points de $\P^2$; on explique comment obtenir les courbes de ramification  ``en haut'' et ``en bas'' ainsi que le moyen de r\'ecup\'erer les $28$ bitangentes \`a partir des sept points de d\'epart (voir section \ref{hyp}). Puis on 
d\'ecrit l'\'eclat\'e du plan projectif le long de sept points par le syst\`eme de cubiques comme hypersurface de bidegr\'e $(2,1)$ de la vari\'et\'e d'incidence ``points-droites'' de $\P^2$. Ceci permet de calculer les images directes des fibr\'es structuraux $O_{\P^2}(n)$, que l'on notera $E_n$, ainsi que leur r\'esolutions et quelques suites exactes dignes d'int\'er\^et (voir la section \ref{En}). Dans la section \ref{tid}, avant d'\'enoncer le th\'eor\`eme concernant les droites de saut des fibr\'es $E_n$, on rappelle le calcul des classes de Chern de toutes les images directes (c'est le thm 5 de \cite{S}) des faisceaux inversibles sur $X_2$, pr\'ealablement d\'efinis \`a partir du groupe de Picard de cette surface.
D\'eterminer les faisceaux inversibles qui ont m\^eme image (au d\'ecalage pr\`es) ne semble pas tr\`es simple. La proposition que Schwarzenberger donne para\^{\i}t inexacte. Nous y revenons dans la derni\`ere section \ref{lesmemes}. 

\subsection{Retour sur l'involution de Geiser}\label{hyp}
\'Etant donn\'es sept points  en position assez g\'en\'erale dans le plan projectif, on consid\`ere l'application rationnelle qui associe, \`a un huiti\`eme point du plan, l'unique pinceau de cubiques passant par les huit points. Si on note $Z=\{ x_1,\cdots,x_7\}$ et $(\Delta_0, \Delta_1, \Delta_2) $ une base du  r\'eseau de cubiques passant par les sept points $x_i$, l'application est tout simplement 
$$ \P^2 ---\rightarrow \P^2, x\mapsto ((\Delta_0(x), \Delta_1(x), \Delta_2(x))$$
Deux cubiques se coupant en $9$ points $\{x,y, x_1,\cdots,x_7\}$, les points $x$ et $y$ ont m\^eme image  et l'application rationnelle  est  de degr\'e deux (l'involution de Geiser, proprement dite, \'echange $x$ et $y$). D\'ecrivons pr\'ecis\'ement l'\'eclatement du plan le long de ce groupe de points $Z$.

\medskip

  On pose
$\P^2=\P V$ et
$W=H^0J_{Z}(3)$. \'Ecrivons une  r\'esolution minimale de cet id\'eal
$$ \begin{CD}
0 @>>> O_{\P^{2}}(-1)\oplus O_{\P^{2}}(-2) @>M>>W\otimes O_{\P^{2}}
@>(\Delta_0,\Delta_1,\Delta_2)>> J_{Z}(3)
@>>> 0
\end{CD}
$$
o\`u, quitte \`a fixer des bases, $M=  \left (
\begin{array}{cc}
            X_0 &C_0\\
              X_1 &C_1\\
               X_2& C_2
\end{array}
  \right )$. On a donc le plongement suivant
$$ \P J_{Z} \hookrightarrow \P W \times \P V =\P^2\times \P^2 $$
L'\'eclatement de $\P^2$ le long de $Z$ est d\'efini dans ce produit par
les \'equations  $$\sum \Delta_iX_i=0 \,\, \mathrm{et} \,\, \sum \Delta_iC_i=0$$
L'image d'un point $x\in \P V$ est le pinceau de cubiques passant par  ce point $x$ et par $Z$. La courbe des cubiques singuli\`eres du r\'eseau $<\Delta_0,\Delta_1,\Delta_2>$ vue dans l'espace des cubiques est une douzique. Cette douzique poss\`ede $28$ points doubles correspondant aux $7$ cubiques doubles en un des $x_i$ mais aussi aux $21$ cubiques d\'ecompos\'ees en une droite passant par deux des sept points et une conique passant par les cinq autres. Cette cubique, doublement singuli\`ere  donne un point double. Chaque point double de la douzique correspond \`a  une bitangente de la courbe duale dans le $\P^2$ des pinceaux de cubiques. On trouve ainsi notre quartique de ramification et ses $28$ bitangentes, dont les sept de d\'epart forment un syst\`eme d'Aronhold, puisqu'elles engendrent les $21$ autres (voir \cite{D} chap.6, pour une d\'efinition plus pr\'ecise d'un syst\`eme d'Aronhold).

\smallskip

On notera $l_i, i=1,\cdots, 7$ les bitangentes correspondantes aux points $x_i, i=1,\cdots, 7$ et $l_{ij}$ celles associ\'ees au choix de deux points $x_i$ et $x_j$ parmi les sept.

\smallskip

La fibre
$q^{-1}(\underline{\alpha})$ est d\'efinie dans $\P V$  par les \'equations 
$$\sum \alpha_iX_i=0 \,\, \mathrm{et} \,\, \sum \alpha_iC_i=0$$ La ramification ``en bas'' (dans $\P W$) correspond aux points 
$\underline{\alpha}$ pour lesquels : 
\begin{center}
$\sum \alpha_iX_i=0$ est
tangente \`a la conique d'\'equation
$\sum \alpha_iC_i=0$. 
\end{center}
 Soit $A=(l_{ij})$ la matrice associ\'ee \`a la conique
$\sum \alpha_iC_i=\sum l_{ij}(\underline{\alpha})X_iX_j$ et $\tilde{A}$ sa
matrice transpos\'ee des
cofacteurs. La droite $\sum \alpha_iX_i=0$ touche la conique d'\'equation
$\sum \alpha_iC_i=0$ si et seulement si $\underline{\alpha}=(\alpha_0,\alpha_1,\alpha_2)$  appartient  \`a la conique duale i.e. $<\tilde{A}\underline{\alpha},\underline{\alpha}>=0$. C'est
l'\'equation de la quartique de
ramification. Cette quartique est lisse en g\'en\'eral (i.e. quand les sept
points sont en position
g\'en\'erale). Bateman (\cite{B}, page 360)  pr\'ecise qu'elle a un point double lorsque trois points sont align\'es ou six sur une conique, qu'elle est form\'ee de deux coniques lorsque six points sont les sommets de quatre droites etc.

\smallskip

Dans $\P V$, ou encore ``en haut'',  la ramification est d\'efinie 
ensemblistement par
$$\{x\in \P^2 \mid H^0(J_{Z}\otimes {\frak m}_x^{2}(3))\neq 0\}$$ C'est la
sextique d'\'equation
$$   {\rm det} \left (
\begin{array}{ccc}
            \frac{\partial \Delta_0}{\partial X_0} &\frac{\partial
\Delta_1}{\partial X_0} & \frac{\partial
\Delta_2}{\partial X_0}\\
              \frac{\partial \Delta_0}{\partial X_1} &\frac{\partial
\Delta_1}{\partial X_1}&\frac{\partial
\Delta_2}{\partial X_1}\\
       \frac{\partial \Delta_0}{\partial X_2}         &\frac{\partial
\Delta_1}{\partial X_2}&
\frac{\partial \Delta_2}{\partial X_2}
\end{array}
  \right )=0  $$ 
Elle est, par d\'efinition, birationnelle \`a la ramification ``en bas'', c'est \`a dire \`a la quartique; elle est donc irr\'eductible et  poss\`ede  $7$
points doubles d'apr\`es Hurwitz.

\subsection{Premiers calculs d'images directes}\label{En}

L'\'eclat\'e $ X=\P J_{Z} $ est d\'efini dans $\P V\times \P W =\P^2\times \P^2$ par deux \'equations 
$\sum X_i\Delta_i=0 $ et $\sum C_i\Delta_i=0$.
La premi\`ere des deux, identifant $W$ \`a $V^{*}$ r\'ealise une dualit\'e. On peut alors consid\'erer que $ \P J_{Z} $ est une hypersurface de bidegr\'ee $(2,1)$ de la vari\'et\'e d'incidence point-droite de $\P V$ o\`u les droites de $\P V$ sont les cubiques passant par les sept points base du d\'epart.

\smallskip

On note $\F \subset \P V \times \P W$ la
vari\'et\'e d'incidence
``points-droites'' de $\P^2$,  $p$ et $q$ les projections,  $\overline{p}$ et $\overline{q}$ leurs restrictions  \`a la sous-vari\'et\'e $\P J_{Z}$.
 $$ \begin{CD}
\P^{2\vee}=\P W  @<\overline{q}<< \P J_{Z} \subset \F  @>q>> \P^{2\vee}=\P W \\
 @. @V\overline{p}VpV \\
@. \P^2=\P V
\end{CD}
$$
Gr\^ace \`a la suite exacte ci-dessous, exprimant que $\P J_{Z}$ est une hypersurface de $\F$ 
$$ \begin{CD}
0 @>>> p^{*}O_{\P^{2}}(-2)\otimes q^{*}O_{\P^{2}}(-1) @>\sum C_i\Delta_i>> O_{\F}
@>>> O_{\P J_{Z}}
@>>>0
\end{CD}
$$
on r\'esout les fibr\'es images directes $\overline{q}_{*}\overline{p}^{*}O_{\P^2}(n)$. Ces images directes  sont de rang deux car elles proviennent de fibr\'es inversibles sur la surface $ \P J_{Z}$ qui est un rev\^etement double de notre plan. 
\begin{prop}\label{resol}
On pose $E_n:=\overline{q}_{*}\overline{p}^{*}O_{\P^2}(n)$. Pour  $n\ge 2$ on a des r\'esolutions
$$ \begin{CD}
0 @>>> [S^{n-2}\Omega^{\vee}(-1)](-1) @>>> S^n\Omega^{\vee}(-1)
@>>> E_n
@>>>0
\end{CD}
$$
De plus,  
$E_0=O_{\P^2}\oplus O_{\P^2}(-2)$, $E_1=\Omega^{\vee}(-1) $ et $E_{-n}=(E_n)^{\vee}(-2)=E_n(-3n)$.\\
Lorsque  $n>0$, on a $H^0(E_n(1-n))=0$, en particulier les fibr\'es $E_n$  sont stables.
\end{prop}
\begin{proof}
La vari\'et\'e d'incidence est une sous-vari\'et\'e du produit $\P^2\times \P^{2\vee}$ d\'efinie par une \'equation de bidegr\'e $(1,1)$ plus suitexactement, 
$$ \begin{CD}
0 @>>> O_{\P^2\times\P^{2\vee}}(-1,-1) @>\sum X_i\Delta_i>> O_{\P^2\times\P^{2\vee}}
@>>> O_{\F}
@>>>0
\end{CD}
$$ 
On calcule l'image directe de $p^{*}O_{\P^2}(n)$ par le morphisme $q$
{\tiny
$$ \begin{CD}
0 @>>> H^0(O_{\P^2}(n-1))\otimes O_{\P^{2\vee}}(-1) @>>> H^0(O_{\P^2}(n-1))\otimes O_{\P^{2\vee}}(-1)
@>>> q_{*}p^{*}O_{\P^2}(n)
@>>>0
\end{CD}
$$ }
Ce qui prouve que $q_{*}p^{*}O_{\P^2}(n)=S^n\Omega^{\vee}(-1)$. C'est ainsi que s'obtient la r\'esolution type. 

\medskip

Le calcul via l'espace produit  $\P^2\times\P^{2\vee}$ du faisceau $R^1q_{*}p^{*}O_{\P^2}(-n)$ permet de montrer 
$E_{-n}=(E_n)^{\vee}(-2)$. Quant au dernier point, on tensorise la r\'esolution de $E_n$ par $O_{\P^2}(1-k)$ pour $k\ge 1$  
$$ \begin{CD}
0 @>>> [S^{n-2}\Omega^{\vee}(-1)](-k) @>>> [S^n\Omega^{\vee}(-1)](1-k)
@>>> E_n(1-k)
@>>>0
\end{CD}
$$
et, en prenant la longue suite exacte de cohomologie, il vient 
$$ \begin{CD}
0 @>>> H^0(E_n(1-k))@>>> H^1([S^{n-2}\Omega^{\vee}(-1)](-k)) @>>> H^1([S^{n}\Omega^{\vee}(-1)](1-k))
\end{CD}
$$
En outre, on v\'erifie facilement que 
$$ \begin{CD}
0 @>>> H^1([S^{n-2}\Omega^{\vee}(-1)](-k)) @>>> S^{n-3}V\otimes S^{k-2}V^{*}
@>>> S^{n-2}V\otimes S^{k-3}V^{*}
\end{CD}
$$
et que cette  derni\`ere fl\`eche est un isomomorphisme quand $k=n$, c'est \`a dire  $$H^1([S^{n-2}\Omega^{\vee}(-1)](-n))=0$$ D'o\`u $H^0(E_n(1-n))=0$
\end{proof}

--- Pour $n=1$ on retrouve $\Omega^{\vee}(-1)$.

\smallskip

--- Pour $n=2$ on trouve le fibr\'e 
$$ \begin{CD}
0 @>>> O_{\P^2}(-1) @>>> S^2\Omega^{\vee}(-1)
@>>> \overline{q}_{*}\overline{p}^{*}O_{\P^2}(2)
@>>>0
\end{CD}
$$ dont la courbe de saut est une courbe de degr\'e $6$ qui est naturellement isomorphe \`a la sextique  d\'efinie ci-avant (la ramification ``en haut''). 
Une r\'esolution minimale de ce fibr\'e est la suivante 
$$ \begin{CD}
0 @>>> O_{\P^2}^4(-1) @>>> O_{\P^2}^6
@>>> \overline{q}_{*}\overline{p}^{*}O_{\P^2}(2)
@>>>0
\end{CD}
$$
C'est le fibr\'e logarithmique (voir \cite{DK}) associ\'e aux sept droites $l_i, i=1,\cdots, 7$.

\smallskip

--- Pour $n=3$ le fibr\'e $E_3$ a pour classes de Chern $(-1,15)$. Son sch\'ema des droites de saut,  bien que support\'e par les bitangentes \`a la quartique de ramification, est de longueur \'egale \`a  $105$. Je ne sais pas s'il n'a que les $7$ droites $l_i, i=1,\cdots, 7$ comme droites de saut ou si toutes  les $28$ bitangentes sont des  droites de saut (dans le thm. \ref{qpon} on montre qu'elles sont toutes de saut pour $E_n$ \`a partir de $n>3$).

\begin{prop}
Soit $l$ une droite g\'en\'erale de $\P^2$. Son image $ \mathcal{L}_n=\overline{q}_{*}\overline{q}^{*}O_{l}(n)$  est   un faisceau de rang un support\'e par la cubique singuli\`ere
$\overline{\lbrace (\Delta_0(x),\Delta_1(x),\Delta_2(x)),x\in l\rbrace} $. Pour tout $n$ il existe des suites exactes 
$$ \begin{CD}
0 @>>>  E_{n-1} @>>>  E_n
@>>> \mathcal{L}_n
@>>>0
\end{CD}$$
\end{prop}
\begin{proof}
Les suites canoniques 
$$ \begin{CD}
0 @>>> O_{\P^2}(n-1) @>>> O_{\P^2}(n)
@>>> O_{l}(n)
@>>>0
\end{CD}$$
remont\'ees par $\overline{p}$ et descendues par $\overline{q}$ deviennent
$$ \begin{CD}
0 @>>>  E_{n-1} @>>> E_{n}
@>>> \mathcal{L}_n
@>>>0
\end{CD}$$
Explicitons, pour $n\ge 1$, le faisceau  $\mathcal{L}_n$. Soit $x$ tel que $x^{\vee}=l$. Remontons la suite exacte 
$$ \begin{CD}
0 @>>> O_{\P^2}(n-1) @>>> O_{\P^2}(n)
@>>> O_{l}(n)
@>>>0
\end{CD}$$
\`a la vari\'et\'e d'incidence et descendons la 
$$ \begin{CD}
0 @>>> S^{n-1}\Omega^{\vee}(-1)@>{l=x^{\vee}}>> S^n\Omega^{\vee}(-1)
@>>> {\frak m}_x^{n}(n)
@>>>0
\end{CD}
$$ 
Maintenant,   utilisons la r\'esolution de $\P J_Z$ pour  trouver 
$$ \begin{CD}
0 @>>> {\frak m}_x^{n-2}(n-3) @>>>  {\frak m}_x^{n}(n) 
@>>> \mathcal{L}_n
@>>>0
\end{CD}$$
 qui d\'efinit concr\`etement le faisceau $\mathcal{L}_n$.
\end{proof}

\subsection{Tous les faisceaux inversibles et leurs images}\label{tid}

Plus g\'en\'eralement nous pouvons, suivant Schwarzenberger, d\'ecrire les faisceaux inversibles sur 
$X_2=\P J_{Z}$ en rappelant que  $\mathrm{Pic}(\P J_{Z})=\Z.l_1+\cdots+\Z.l_7+\Z.h$ o\`u 
les $l_i$ sont les droites (diviseurs exceptionnels) images inverses des sept points $x_i$,  et $h=p^{*}H$ avec $H$ qui est une section hyperplane g\'en\'erale de $\P V$. Les relations entre les g\'en\'erateurs du groupe de Picard sont  
$h^2=1, h.l_i=0, l_il_j=0$ et $l_i^2=-1$. Le diviseur $nh+\sum t_il_i$ induit un fibr\'e en droites not\'e $L^{n,t_i}$. Ses classes de Chern valent
\begin{prop}(\cite{S}, thm. 6)\\
$c_1(\overline{q}_{*}L^{n,t_i})=\sum t_i+3n-2$ \\
$c_2(\overline{q}_{*}L^{n,t_i})=4n^2-3n+(3n-1)(\sum t_i)+(\sum t_i^2)+(\sum_{i<j}t_it_j)$
\end{prop}
Comme nous l'avons rappel\'e dans la partie consacr\'ee \`a l'involution de Geiser, les sept points $x_i$ ont pour images les droites bitangentes $l_i$ tandis que chaque couple de points $(x_i,x_j)$ d\'efinit une unique droite $d_{ij}$ dont l'image est aussi une droite, not\'ee $l_{ij}$ qui est une des $21$ bitangentes restantes. Ces $28$ droites apparaissent dans les sch\'emas de droites de saut comme d\'ecrit dans le th\'eor\`eme ci-dessous.
\begin{theo} \label{qpon}Soit $n\ge 2$, les $28$ bis\'ecantes sautent de la mani\`ere suivante, \\
1) La droite $l_i$ ``image'' du point $x_i$ est une droite de saut d'ordre $\lceil \frac{3n-2}{2} \rceil$ pour le fibr\'e $E_n$.\\
2) La droite $l_{i,j}$ ``image'' de la droite joignant $x_i$ et $x_j$ est une droite de saut d'ordre 
$\lceil \frac{n-2}{2} \rceil$ pour le fibr\'e $E_n$.
\end{theo}
\begin{rem}\label{equivalence}
Lorsque  $n\ge 4$, les $28$ bis\'ecantes sont des droites sauteuses des fibr\'es $E_n$, et elles sont les seules lorsque $n$ est impair.
\end{rem}
La d\'emonstration de ce th\'eor\`eme repose sur l'existence de suites exactes, que l'on nommera  \textit{suites de liaison}, provenant des suites exactes canoniques sur $\P^2$ liant les faisceaux ${\frak m}_{x_i}(n), J_{x_i,x_j}(n)$ et $O_{\P^2}(n)$ entre-eux. En passant on montrera que la droite $l_{i,j}$ image de la droite joignant $x_i$ et $x_j$ est une droite de saut d'ordre 
$\lceil \frac{n}{2} \rceil$ pour le fibr\'e $\overline{q}_{*}\overline{p}^{*}J_{x_i,x_j}(n)$ et d'ordre de saut 
$\lceil \frac{n-1}{2} \rceil$ pour les fibr\'es $\overline{q}_{*}\overline{p}^{*}{\frak m}_{x_i}(n)$ et 
$\overline{q}_{*}\overline{p}^{*}{\frak m}_{x_j}(n)$. Avec les m\^emes techniques on peut montrer que la droite 
$l_i$ est une droite de saut d'ordre $\lceil \frac{3n-3}{2} \rceil$ pour le fibr\'e $\overline{q}_{*}\overline{p}^{*}{\frak m}_{x_j}(n)$.

\begin{proof}
La suite canonique
$$ \begin{CD}
0 @>>> {\frak m}_{x_i}(n)  @>>> O_{\P^2}(n)
@>>> O_{x_i}
@>>>0
\end{CD}
$$
aura pour image directe la suite 
$$ \begin{CD}
0 @>>> \overline{q}_{*}\overline{p}^{*}{\frak m}_{x_i}(n)  @>>> E_n
@>>> \overline{q}_{*}\overline{p}^{*}O_{x_i}
@>>>0
\end{CD}
$$
Le dernier faisceau est support\'e par $\overline{q}\overline{p}^{-1}(x_i)=l_i$; ce dont on s'assure en v\'erifiant que 
$c_1(E_n)=c_1(\overline{q}_{*}\overline{p}^{*}{\frak m}_{x_i}(n))+1$. Enfin le calcul de la seconde classe de Chern donne le d\'ecalage, i.e $ \overline{q}_{*}\overline{p}^{*}O_{x_i}=O_{l_i}$. Par cons\'equent  $E_{n\mid l_i}=O_{l_i}\oplus O_{l_i}(3n-2)$, ce qui prouve le premier point.

\smallskip

En ce qui concerne le second point on reprend les m\^emes id\'ees et techniques : suites de liaison et calcul des classes de Chern afin de conna\^itre le degr\'e du support et le d\'ecalage. Ainsi les suites 
$$ \begin{CD}
0 @>>> O_{\P^2}(n-1)  @>>> J_{x_i,x_j}(n)
@>>> O_{d_{i,j}}
@>>>0
\end{CD}
$$
auront  pour images directes les suites 
$$ \begin{CD}
0 @>>> E_{n-1}  @>>> \overline{q}_{*}\overline{p}^{*}J_{x_i,x_j}(n)
@>>> O_{l_{i,j}}(n-2)
@>>>0
\end{CD}
$$
Ce qui implique $[\overline{q}_{*}\overline{p}^{*}J_{x_i,x_j}(n)]_{\mid l_{i,j}}=O_{l_{i,j}}(n-2)\oplus O_{l_{i,j}}(2n-2)$.

\smallskip

Les suites, 
$$ \begin{CD}
0 @>>>  J_{x_i,x_j}(n) @>>> {\frak m}_{x_i}(n)
@>>> O_{x_j}
@>>>0
\end{CD}
$$
auront  pour images directes les suites 
$$ \begin{CD}
0 @>>> \overline{q}_{*}\overline{p}^{*}J_{x_i,x_j}(n)  @>>> \overline{q}_{*}\overline{p}^{*}{\frak m}_{x_i}(n)
@>>> O_{l_j}
@>>>0
\end{CD}
$$
Apr\`es tensorisation des faisceaux de  cette suite exacte par $O_{l_{ij}}$ on observe que 
$$ \overline{q}_{*}\overline{p}^{*}{\frak m}_{x_i}(n)\otimes O_{l_{ij}}= O_{l_{ij}}(n-2)\oplus O_{l_{ij}}(2n-1)\,\,
\mathrm{ou} \,\, O_{l_{ij}}(n-1)\oplus O_{l_{ij}}(2n-2)$$
Enfin pour conclure on restreint la suite exacte 
$$ \begin{CD}
0 @>>> \overline{q}_{*}\overline{p}^{*}{\frak m}_{x_i}(n)  @>>> E_n 
@>>> O_{l_i}
@>>>0
\end{CD}
$$
\`a la droite $l_{i,j}$.
\end{proof}
Ces d\'evissages via les suites de liaison permettent d'\'etudier les fibr\'es images directes pas \`a pas, la connaissance de l'un fournissant des renseignements sur le suivant. Ainsi pour tout $n$, le fibr\'e $E_n$ vient avec sa cohorte de sept compagnons, $E_n^{k}:=q_{*}p^{*}J_{x_1,\cdots,x_k}(n)$ v\'erifiant 
$$ \begin{CD}
0 @>>> E_n^{k} @>>> E_n 
@>>> \oplus_{i=1}^{i=k} O_{l_i}
@>>>0
\end{CD}
$$ 
D\'etaillons un ou deux exemples.

\smallskip

\textbf{Exemple 1)} Pour $E_2$ qui est le fibr\'e logarithmique (voir \cite{DK}) associ\'e aux sept droites $l_i$ on trouve $E_2^{1}\in M(-1,4)$ qui est le fibr\'e logarithmique associ\'e aux six droites restantes, $E_2^{2}\in M(0,2)$ dont les droites de saut sont les droites tangentes \`a la conique d\'efinie par les cinq droites restantes, 
$E_2^{3}=\Omega^{\vee}(-1)$, $E_2^{4}=O_{\P^2}^2$, $E_2^{5}=O_{\P^2}\oplus O_{\P^2}(-1)$, 
$E_2^{6}=O_{\P^2}^2(-1)$ et enfin $E_2^{7}=\Omega$.
Ainsi les r\'esolutions minimales se d\'eduisent les unes des autres et pour $k\le 4$ on peut \'ecrire
$$ \begin{CD}
0 @>>> O_{\P^2}^{4-k}(-1) @>>> O_{\P^2}^{6-k}
@>>> E_2^k
@>>>0
\end{CD}
$$
\textbf{Exemple 2)} Comme deuxi\`eme exemple regardons le cas de $E_3$. Il admet une  r\'esolution minimale du type
 $$ \begin{CD}
0 @>>> O_{\P^2}^6(-1) @>>> O_{\P^2}^7\oplus O_{\P^2}(1)
@>>> E_3
@>>>0
\end{CD}
$$
Consid\'erons l'unique  section non nulle de $E_3(-1)$. Elle donne 
 $$ \begin{CD}
0 @>>> O_{\P^2} @>>> E_3(-1)
@>>> J_{\Gamma}(5)
@>>>0
\end{CD}
$$
o\`u $\Gamma$ est un groupe de points de longueur \'egale \`a $21$.
En la sym\'etrisant il vient
 $$ \begin{CD}
0 @>>> E_3(-1) @>>> (S^2E_3)(-2)
@>>> J_{\Gamma}^2(10)
@>>>0
\end{CD}
$$
\`A la section non nulle de  $H^0((S^2E_3)(-5))\neq 0$, donn\'ee par la quartique de ramification, correspond une courbe de degr\'e $7$ qui poss\`ede les $21$ points de $\Gamma$ comme points doubles. Ce ne peut-\^etre qu'une configuration de sept droites. Celles-ci \'etant \'evidemment les droites $l_i$. Alors 
$E_3^{1}\in M(0,11)$ et $h^0(E_3^{1}(-1))=1$ et les sommets des six droites restantes sont les z\'eros de son unique section. De m\^eme $E_3^{2}\in M(-1,8)$, $h^0(E_3^{2}(-1))=1$ et les sommets des cinq droites restantes sont  les z\'eros de son unique section, etc. Le dernier fibr\'e \'etant $\overline{q}_{*}\overline{p}^{*}J_Z(3)=O_{\P^2}(-1)\oplus O_{\P^2}(1)$.
Ainsi les r\'esolutions se d\'eduisent les unes des autres et pour $k\le 6$ on peut \'ecrire
$$ \begin{CD}
0 @>>> O_{\P^2}^{6-k}(-1) @>>> O_{\P^2}^{7-k}\oplus O_{\P^2}(1)
@>>> E_3^k
@>>>0
\end{CD}
$$
\`A chaque pas la section de $S^2E_3^k(2-c_1)$ provenant de la quartique de ramification se traduit par l'existence d'une configuration de droites. Il faut remarquer que la quartique n'est parfaitement d\'efinie, ou encore  $h^0(S^2E_3^k(2-c_1))=1$ que pour $E_3$, apr\`es il n'y a plus assez de droites. 

\smallskip

Disons aussi un mot de la th\'eta-caract\'eristique associ\'ee dans le cas de $E_3$. On a une suite exacte canonique provenant du rev\^etement 
$$\begin{CD}
0 @>>>E_3(-2) @>>>E_3
@>>> \theta
@>>>0
\end{CD}
$$
avec $\theta(-8)=\theta^{*}\otimes O_{C_{4}}(1)$. Soit encore 
$\theta^2=O_{C_{4}}(9)$. L'unique section de $E_3(-1)$ correspond \`a une section de $\theta(-1)$, ou encore \`a une racine carr\'ee de $O_{C_{4}}(7)$ dont les z\'eros sont les quatorze points de contact des $7$ droites bitangentes.

\begin{rem}
Tandis que tous les fibr\'es de Schwarzenberger poss\`edaient une r\'esolution minimale par une matrice de formes lin\'eaires, il n'y a pas, pour les fibr\'es $\overline{q}_{*}L^{n,t_i}$, ni m\^eme pour les fibr\'es $E_n$ de r\'esolution type. 
Ainsi on semble avoir 
$$ \begin{CD}
0 @>>> O_{\P^2}^{7}(-1) @>>> O_{\P^2}^{3}(1)\oplus O_{\P^2}^6
@>>> E_4
@>>>0
\end{CD}
$$
tandis que 
$$ \begin{CD}
0 @>>> O_{\P^2}^{7}(-1) @>>> O_{\P^2}^{6}(1)\oplus O_{\P^2}^3
@>>> E_5
@>>>0
\end{CD}
$$
et encore 
$$ \begin{CD}
0 @>>> O_{\P^2}^{6}(-1)\oplus O_{\P^2}^2 @>>> O_{\P^2}^{10}(1)
@>>> E_6
@>>>0
\end{CD}
$$
\end{rem}

\subsection{Ceux qui ont la m\^eme image}
\label{lesmemes}
Soucieux d'\'etudier les restrictions des fibr\'es $\overline{q}_{*}L^{n,t_i}$ aux droites du plan, Schwarzenberger affirme \`a juste titre qu'il est n\'ec\'essaire, au pr\'ealable, de reconna\^{\i}tre les fibr\'es en droites donnant le m\^eme fibr\'e de rang deux. Il donne  sur ce point un \'enonc\'e (prop 9, page 633) qui semble inexact puisqu'il impliquerait, par exemple, $\overline{q}_{*}L^{1,0}=\overline{q}_{*}L^{0,-3}$; or le premier est le fibr\'e  tangent tandis que le deuxi\`eme a pour classes de Chern, $c_1=-1, c_2=631$ apr\`es normalisation.

\smallskip

On se contentera ici, dans la  proposition \ref{meme}, d'identifier les faisceaux inversibles avec le m\^eme poids $t$ sur chaque $x_i$ ayant au d\'ecalage pr\`es la m\^eme image directe sur le plan projectif. Pour cela, nous cherchons dans un premier temps \`a 
\'ecrire $\overline{q}^{*}O_{\P^2}(m)$ dans la base de $Pic(X_2)$ c'est \`a dire \`a trouver $t_i$ et $n$ tels que $\overline{q}^{*}O_{\P^2}(m)=L^{n,t_i}$; ceci permettra connaissant $L$ et son image $E$ d'en d\'eduire le faisceau inversible sur $X_2$ dont l'image est un d\'ecal\'e $E(m)$ de $E$ (proposition \ref{meme} ci-dessous) 
\begin{lem} Pour tout $m\in \Z$ on a 
 $\overline{q}^{*}O_{\P^2}(m)=L^{3m,-m}$
\end{lem}
\begin{proof}
Comme $\overline{q}_{*}\overline{q}^{*}O_{\P^2}(m)=O_{\P^2}(m)\oplus O_{\P^2}(m-2)$ on cherche \`a r\'esoudre les \'equations 
$$\left \lbrace
\begin{array}{l}
2m-2=3n-2+\sum t_i\\
m^2-2m=4n^2-3n+(3n-1)(\sum t_i)+(\sum t_i^2)+(\sum_{i<j}t_it_j)
\end{array}
\right.
$$
Aucun des sept points n'\'etant privil\'egi\'e tous les $t_i$ sont \'egaux et on note $t=t_i$ la valeur commune. Alors le syst\`eme d'\'equation ci-dessus \'equivaut \`a 
$$\left \lbrace
\begin{array}{l}
2m=3n+7t\\
m^2=4n^2+21nt+28t^2
\end{array}
\right.
$$
En faisant la diff\'erence de $4$ fois la deuxi\`eme  et du carr\'e de la premi\`ere on se ram\`ene \`a 
$$n^2+6nt+9t^2=(n+3t)^2=0$$ 
En remontant il vient $t=-m$ et $n=3m$.
\end{proof}
\begin{rem}
Avec les m\^emes techniques on montre qu'il n'existe pas de fibr\'e $L^{n,t_i}$ avec tous les $t_i$ \'egaux tel que les classes de Chern du normalis\'e soit $c_1=-1$ et $c_2=8$. C'est dommage car un tel fibr\'e  poss\`edant un sch\'ema de droites de saut de longueur \'egale \`a $28$, aurait pu  avoir pour droites de saut exactement (sans multiplicit\'e) les  $28$ bitangentes offrant ainsi  une m\'ethode pour reconstruire la quartique \`a partir des bitangentes. 
\end{rem}

Cette lemme et la formule de projection ont pour cons\'equence imm\'ediate
\begin{prop}\label{meme}
$\overline{q}_{*}L^{n+3m,t_i-m}=[\overline{q}_{*}L^{n,t_i}](m)$
\end{prop}
Par exemple
$[\overline{q}_{*}\overline{p}^{*}{\frak m}_{x_1}(2)](1)=\overline{q}_{*}\overline{p}^{*}[J_{Z\setminus \lbrace x_1\rbrace}\otimes {\frak m}_{x_1}^2]$ et plus g\'en\'eralement 
$$E_{n-3}(1)=\overline{q}_{*}\overline{p}^{*}J_{Z}(n)$$

Vall\`es Jean\\
Laboratoire de Math\'ematiques appliqu\'ees\\
de Pau et de Pays de l'Adour,\\
 Avenue de l'Universit\'e\\
 64000 Pau(France)\\
 email:  jean.valles@univ-pau.fr \\

\end{document}